\documentclass{article}
\usepackage{amssymb,amsmath,epsfig,graphics}

\newcounter{sscounter}[section]
\renewcommand{\thesscounter}{\thesection.\arabic{sscounter}}
\newcommand{\ssnnl}[1]{\noindent
\refstepcounter{sscounter}\bf\thesscounter. #1\rm}

\newcommand{\CC}{{\mathbb C }}
\newcommand{\KK}{{\mathbb K }}

\newcommand{\QQ}{{\mathbb Q }}
\newcommand{\RR}{{\mathbb R }}
\newcommand{\TT}{{\mathbb T }}
\newcommand{\ZZ}{{\mathbb Z }}

\newcommand{\cA}{\mathcal{A}}

\newcommand{\cP}{\mathcal{P}}
\newcommand{\cQ}{\mathcal{Q}}

\newcommand{\cZ}{\mathcal{Z}}

\newcommand{\va}{\mathbf{a}}
\newcommand{\vb}{\mathbf{b}}

\newcommand{\vt}{\mathbf{t}}
\newcommand{\vx}{\mathbf{x}}
\newcommand{\vz}{\mathbf{z}}
\newcommand{\vu}{\mathbf{u}}
\newcommand{\vv}{\mathbf{v}}
\newcommand{\vw}{\mathbf{w}}

\newcommand{\nv}{\mathbf{0}}
\newcommand{\rank}{\mathrm{rank}}

\newcommand{\ga}{\alpha}
\newcommand{\gb}{\beta}
\newcommand{\gc}{\gamma}
\newcommand{\gf}{\varphi}
\newcommand{\gl}{\lambda}
\newcommand{\gs}{\sigma}

\newcommand{\gz}{\zeta}
\newcommand{\ep}{\varepsilon}

\newcommand{\gB}{\mathfrak{B}}
\newcommand{\gW}{\mathfrak{W}}

\newcommand{\edot}{\wedge}

\newcommand{\modquot}[2]{\mbox{\raisebox{.2ex}{$#1$}\hspace{-.3em}/ \hspace{-.6em} \raisebox{-.2ex}{$#2$}}}

\begin{document}

\title{Computation of Principal $\cA$-determinants through Dimer Dynamics}

\author{Jan Stienstra\\
\small Mathematisch Instituut, Universiteit Utrecht, the Netherlands\\ 
\small e-mail: {J.Stienstra}{`at'}{uu.nl} \normalsize}

\date{}

\maketitle

\abstract{$\cA$ is a set of $N$ vectors in $\ZZ^{N-2}$
situated in a hyperplane not through $0$ and spanning $\ZZ^{N-2}$ over $\ZZ$.
Gulotta's algorithm \cite{G} constructs from $\cA$ a dimer model. 
A theorem  in \cite{S2} states that the principal $\cA$-determinant 
equals the determinant of (a suitable form of) the Kasteleyn matrix of that dimer model. In the present note we translate Gulotta's pictorial description of the algorithm into matrix operations. As a result one obtains an algorithm for computing
the principal $\cA$-determinant, which is much faster than the algorithm in \cite{S1}.}

\section{Introduction}\label{intro}

$\cA\,=\,\{\va_1,\ldots,\va_N\}$ is a set of vectors
in $\ZZ^{N-2}$ situated in a hyperplane not through $0$ and spanning $\ZZ^{N-2}$ over $\ZZ$.
The principal $\cA$-determinant, defined in \cite{GKZ}, describes the singularities of Gelfand-Kapranov-Zelevinsky's $\cA$-hypergeometric system of partial differential equations \cite{gkz1}. It also describes for which $N$-tuples of coefficients
$u_1,\ldots,u_N\in\CC$ the Laurent polynomial
$\sum_{j=1}^N u_j\vx^{\va_j}$ in $N-2$ variables is singular (see \cite{GKZ} for details). It
is a polynomial with integer coefficients in the variables $u_1,\ldots,u_N$. 
The restriction $\mathrm{rank }\cA=\sharp\cA-2$ means that the corresponding
hypergeometric functions are essentially functions in two variables and that
in the Laurent polynomial the number of terms exceeds the number of variables by $2$. Even in this case the definition of the principal $\cA$-determinant 
is fairly complicated and only a few of its coefficients could explicitly be
calculated in \cite{GKZ}. In \cite{DS} Dickenstein and Sturmfels re-examined the definition of the principal $\cA$-determinant and related it to Chow forms
which is another important concept from \cite{GKZ}.
In \cite{S2} it was shown how these Chow forms and the principal $\cA$-determinant can be easily computed as the determinant of a suitable version of the Kasteleyn
matrix of a dimer model associated with $\cA$.
When writing \cite{S2} I had only the algorithm in \cite{S1} to construct that dimer model. Later Gulotta gave another algorithm \cite{G} for constructing the appropriate
dimer model. He describes the algorithm as a process that transforms certain doubly periodic configurations of curves in the plane.
In Sections \ref{section:zigzag}, \ref{section:algorithm1}, \ref{section:algorithm2} we give a faithfull reproduction of those configurations of curves
by matrices and of Gulotta's algorithm by row and column operations on these matrices.
In that form Gulotta's algorithm, which is a fast converging iterative process,
is much more efficient than the algorithm in \cite{S1}, which is a search with many trial-and-errors. 
Moreover, unlike for the algorithm in \cite{S1} for Gulotta's algorithm it can be guaranteed that it finds a desired dimer model. On the other hand there are cases in which \cite{S1} yields two different models and \cite{G} gives only one.

In this note we do not need formal definitions of `dimer model'
and `principal $\cA$-determinant'. Dimer models are implicitly present through their patterns of zigzags. This is briefly explained in 
Remark \ref{rem:zigzag dimer}. Principal $\cA$-determinants appear only in 
Section \ref{sec:Adet} in a quotation from \cite{S2}.

In Section \ref{sec:KZ} we recall from \cite{S1} how a pattern of zigzags 
(alias dimer model) is faithfully represented by a matrix 
$\KK_\cZ(\vz,\vu)$, which is in fact a suitable \textit{generalization of the Kasteleyn matrix of the dimer model}.
In Section \ref{sec:Adet} we recall from \cite{S2} the theorem that expresses
the principal $\cA$-determinant as the determinant of $\KK_\cZ(\vz,\vu)$.
Sections \ref{sec:KZ} and \ref{sec:Adet} and  can be read immediately after Definitions \ref{def:zigzag pattern} and \ref{def:good}.

\section{Patterns of zigzags on the torus}
\label{section:zigzag}

\ssnnl{}\label{zigzag pattern}
This note is about patterns of zigzags on the torus $\TT=\modquot{\RR^2}{\ZZ^2}$.
Here \textit{zigzag} means the image (modulo $\ZZ^2$) of an oriented connected curve $C$ in the plane 
(i.e. the image of a continuous map from $\RR$ to $\RR^2$) such that the two coordinate functions restrict to monotone functions on $C$ and such that $\vt+C=C$ for some non-zero vector $\vt\in\ZZ^2$.

We denote by $\ep_1$ (resp. $\ep_2$) the zigzags in $\TT$ coming from the first
(resp. second) coordinate axis of $\RR^2$.
The homology classes of $\ep_1$ and $\ep_2$ form the standard basis for the homology group
$H_1(\TT,\ZZ)$. We denote the intersection number of two elements
$\ga,\,\gb\in H_1(\TT,\ZZ)$ by $\ga\edot\gb$. For two zigzags $Z,\,Z'$ on $\TT$
we write $Z\edot Z'$ for the intersection number of their homology classes.
For two zigzags $Z,\,Z'$ which intersect in only a finite number of points and for which the intersections are transverse, 
each intersection point contributes $+1$ or $-1$ to $Z\edot Z'$ according to the orientation; e.g.
$\ep_1\edot\ep_1\,=\,\ep_2\edot\ep_2\,=\,0$,  
$\:\ep_1\edot\ep_2\,=\,-\ep_2\edot\ep_1\,=\,1$.

Writing $[Z]$ for the homology class of a zigzag $Z$ we have
\begin{equation}\label{eq:Zee}
[Z]\,=\,(Z\edot\ep_2)[\ep_1]\,-\,(Z\edot\ep_1)[\ep_2]\,.
\end{equation}

\

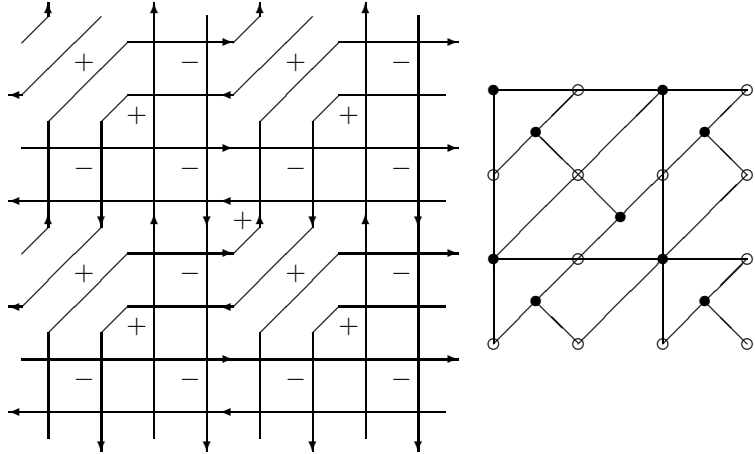
\begin{figure}[ht]
\begin{picture}(160,190)(0,0)
\put(30,0){
\begin{picture}(160,160)(0,0)
\put(0,0){
\begin{picture}(80,80)(0,0)
\put(0,70){\line(1,1){10}}
\put(0,50){\line(1,1){30}}
\put(10,40){\line(1,1){30}}
\put(10,40){\line(0,-1){40}}
\put(40,70){\line(1,0){40}}
\put(30,40){\line(1,1){10}}
\put(30,40){\line(0,-1){40}}
\put(40,50){\line(1,0){40}}
\put(0,30){\line(1,0){80}}
\put(0,10){\line(1,0){80}}
\put(50,0){\line(0,1){80}}
\put(70,0){\line(0,1){80}}
\multiput(10,80)(40,0){2}{\vector(0,1){5}}
\multiput(30,0)(40,0){2}{\vector(0,-1){5}}
\multiput(75,30)(0,40){2}{\vector(1,0){5}}
\multiput(0,10)(0,40){2}{\vector(-1,0){5}}
\end{picture}
}
\put(80,0){
\begin{picture}(80,80)(0,0)
\put(0,70){\line(1,1){10}}
\put(0,50){\line(1,1){30}}
\put(10,40){\line(1,1){30}}
\put(10,40){\line(0,-1){40}}
\put(40,70){\line(1,0){40}}
\put(30,40){\line(1,1){10}}
\put(30,40){\line(0,-1){40}}
\put(40,50){\line(1,0){40}}
\put(0,30){\line(1,0){80}}
\put(0,10){\line(1,0){80}}
\put(50,0){\line(0,1){80}}
\put(70,0){\line(0,1){80}}
\multiput(10,80)(40,0){2}{\vector(0,1){5}}
\multiput(30,0)(40,0){2}{\vector(0,-1){5}}
\multiput(80,30)(0,40){2}{\vector(1,0){5}}
\multiput(0,10)(0,40){2}{\vector(-1,0){5}}
\end{picture}
}
\put(0,80){
\begin{picture}(80,80)(0,0)
\put(0,70){\line(1,1){10}}
\put(0,50){\line(1,1){30}}
\put(10,40){\line(1,1){30}}
\put(10,40){\line(0,-1){40}}
\put(40,70){\line(1,0){40}}
\put(30,40){\line(1,1){10}}
\put(30,40){\line(0,-1){40}}
\put(40,50){\line(1,0){40}}
\put(0,30){\line(1,0){80}}
\put(0,10){\line(1,0){80}}
\put(50,0){\line(0,1){80}}
\put(70,0){\line(0,1){80}}
\multiput(10,80)(40,0){2}{\vector(0,1){5}}
\multiput(30,5)(40,0){2}{\vector(0,-1){5}}
\multiput(75,30)(0,40){2}{\vector(1,0){5}}
\multiput(0,10)(0,40){2}{\vector(-1,0){5}}
\end{picture}
}
\put(80,80){
\begin{picture}(80,80)(0,0)
\put(0,70){\line(1,1){10}}
\put(0,50){\line(1,1){30}}
\put(10,40){\line(1,1){30}}
\put(10,40){\line(0,-1){40}}
\put(40,70){\line(1,0){40}}
\put(30,40){\line(1,1){10}}
\put(30,40){\line(0,-1){40}}
\put(40,50){\line(1,0){40}}
\put(0,30){\line(1,0){80}}
\put(0,10){\line(1,0){80}}
\put(50,0){\line(0,1){80}}
\put(70,0){\line(0,1){80}}
\multiput(10,80)(40,0){2}{\vector(0,1){5}}
\multiput(30,5)(40,0){2}{\vector(0,-1){5}}
\multiput(80,30)(0,40){2}{\vector(1,0){5}}
\multiput(0,10)(0,40){2}{\vector(-1,0){5}}
\end{picture}
}
\put(143,140){$-$}
\put(123,120){$+$}
\put(103,100){$-$}
\put(83,80){$+$}
\put(63,60){$-$}
\put(43,40){$+$}
\put(23,20){$-$}
\put(23,60){$+$}
\put(63,100){$-$}
\put(103,140){$+$}
\put(23,100){$-$}
\put(43,120){$+$}
\put(63,140){$-$}
\put(23,140){$+$}
\put(63,20){$-$}
\put(103,20){$-$}
\put(143,20){$-$}
\put(123,40){$+$}
\put(103,60){$+$}
\put(143,60){$-$}
\put(143,100){$-$}
\end{picture}}

\put(190,20){
\setlength{\unitlength}{.8pt}

\begin{picture}(160,160)(0,0)
\put(143,140){\circle{5}}
\put(123,120){\circle*{5}}
\put(103,100){\circle{5}}
\put(83,80){\circle*{5}}
\put(63,60){\circle{5}}
\put(43,40){\circle*{5}}
\put(23,20){\circle{5}}
\put(23,60){\circle*{5}}
\put(63,100){\circle{5}}
\put(103,140){\circle*{5}}
\put(23,100){\circle{5}}
\put(43,120){\circle*{5}}
\put(63,140){\circle{5}}
\put(23,140){\circle*{5}}
\put(63,20){\circle{5}}
\put(103,20){\circle{5}}
\put(143,20){\circle{5}}
\put(123,40){\circle*{5}}
\put(103,60){\circle*{5}}
\put(143,60){\circle{5}}
\put(143,100){\circle{5}}
\put(23,20){\line(1,1){120}}
\put(63,20){\line(1,1){80}}
\put(103,20){\line(1,1){40}}
\put(23,60){\line(1,1){80}}
\put(23,100){\line(1,1){40}}
\put(23,20){\line(0,1){120}}
\put(103,20){\line(0,1){120}}
\put(23,60){\line(1,0){120}}
\put(23,140){\line(1,0){120}}
\put(143,20){\line(-1,1){20}}
\put(83,80){\line(-1,1){40}}
\put(63,20){\line(-1,1){20}}
\put(143,100){\line(-1,1){20}}
\end{picture}}
\setlength{\unitlength}{1pt}
\end{picture}
\caption{\label{fig:zigzag dimers}
\textit{Zigzag pattern and corresponding dimer model.}}
\end{figure}

\begin{figure}[ht]
\begin{picture}(80,140)(-40,0)
\put(10,20){
\dashbox{1.0}(80,80){
\begin{picture}(80,80)(0,0)
\put(10,80){\vector(0,-1){80}}
\put(50,40){\vector(0,-1){40}}
\put(30,0){\vector(0,1){80}}
\put(70,60){\line(-1,1){10}}
\put(70,0){\line(0,1){40}}
\put(80,70){\vector(-1,1){10}}
\put(40,70){\vector(-1,0){40}}
\put(80,30){\vector(-1,0){80}}
\put(0,50){\line(1,0){40}}
\put(40,50){\line(1,-1){10}}
\put(0,10){\vector(1,0){80}}
\put(70,40){\line(-1,1){30}}
\put(50,80){\vector(1,-1){30}}

\footnotesize{
\put(7,85){$1$}
\put(47,85){$2$}
\put(27,-7){$4$}
\put(67,-7){$5$}
\put(83,69){$5$}
\put(83,29){$6$}
\put(-10,49){$2$}
\put(-10,9){$3$}

\put(55,60){$-$}
\put(55,20){$+$}
\put(15,60){$+$}
\put(15,20){$+$}
\put(0,2){$-$}
\put(0,40){$-$}
\put(0,72){$-$}
\put(37,2){$-$}
\put(37,40){$-$}
\put(75,2){$-$}
\put(75,75){$-$}

}
\end{picture}
}}
\put(95,60){\footnotesize{
$
\begin{array}{rl}
\begin{array}{ll}
&\hspace{-.5em}B_\cZ\\
I_\cZ&\hspace{-.5em}P_\cZ
\end{array}&\hspace{-1.5em}
\left[\!\!\begin{array}{rrrrrr}
1\!\!&1\!\!&0\!\!&-1\!\!&-1\!\!&0\\
0\!\!&1\!\!&1\!\!&0\!\!&-1\!\!&-1
\end{array}\!\!\right]\\ [2ex]
\left[\!\!\begin{array}{l}
110000\\
100010\\
010100\\
000110\\
101000\\
100001\\
001100\\
000101\\
011000\\
010001\\
001010\\
000011
\end{array}\!\!\right]&\hspace{-1.5em}
\left[\!\!\begin{array}{rrrrrr}
1\!\!&0\!\!&0\!\!&0\!\!&1\!\!&0\\
1\!\!&0\!\!&0\!\!&0\!\!&1\!\!&0\\
1\!\!&0\!\!&0\!\!&0\!\!&1\!\!&0\\
1\!\!&0\!\!&0\!\!&0\!\!&1\!\!&0\\
1\!\!&-1\!\!&0\!\!&0\!\!&1\!\!&1\\
1\!\!&-1\!\!&0\!\!&0\!\!&1\!\!&1\\
1\!\!&-1\!\!&0\!\!&0\!\!&1\!\!&1\\
1\!\!&-1\!\!&0\!\!&0\!\!&1\!\!&1\\
1\!\!&0\!\!&0\!\!&-1\!\!&1\!\!&1\\
1\!\!&0\!\!&0\!\!&-1\!\!&1\!\!&1\\
1\!\!&0\!\!&0\!\!&-1\!\!&1\!\!&1\\
1\!\!&0\!\!&0\!\!&-1\!\!&1\!\!&1
\end{array}\!\!\right]
\end{array}
$
}}

\end{picture}
\caption{\label{fig:zigzag matrices}
\textit{Pattern of zigzags and the corresponding matrices.
}}
\end{figure}
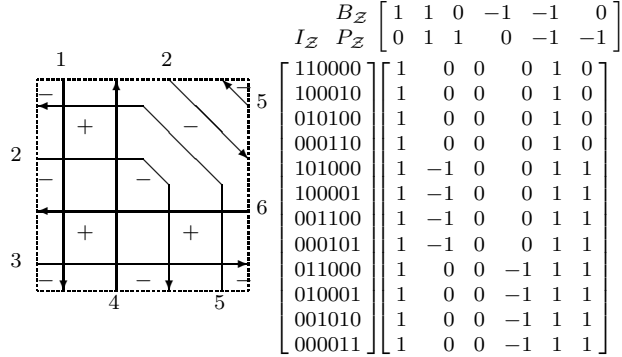

\ssnnl{Definition.}\label{def:zigzag pattern}
In this note \textit{pattern of zigzags} means a finite sequence
$\cZ=(Z_1,\ldots,Z_p)$ of zigzags on $\TT$ which satisfies the following conditions:
\begin{enumerate}
\item
The homology classes $[Z_1],\ldots,[Z_p]$ span
$H_1(\TT,\QQ)$.
Their sum is $0$.
\item
For every $j$ the greatest common divisor of 
$(Z_j\edot\ep_1,Z_j\edot\ep_2)$ is $1$.
\item
Every point of $\TT$ lies on at most two zigzags in $\cZ$.
\item
Every pair of zigzags $Z_i,\,Z_j$ intersects in only a finite number of points and the intersections are transverse. 
\item
The $2$-cells (i.e. connected components) of $\TT\setminus\bigcup_{i=1}^p Z_i$
are divided into three types: those (called $+$-cells) of which the boundary is positively oriented,
those (called $-$-cells) of which the boundary is negatively oriented and those of which the boundary is not oriented. It is required that for every intersection point $x$ 
of zigzags in $\cZ$ one of the four $2$-cells having $x$ in their boundary
is a $+$-cell, one is a $-$-cell
and two have an unoriented boundary; see Figures 
\ref{fig:zigzag dimers}-\ref{fig:start pattern}.
\item
The number of $+$-cells equals the number of $-$-cells.
\end{enumerate}

\

\ssnnl{Definition.}\label{def:good} 
We say that a
pattern of zigzags $\cZ=(Z_1,\ldots,Z_p)$ is \textit{good} if in addition to the above conditions, it also satisfies:
\begin{enumerate}\addtocounter{enumi}{6}
\item
Write $\gz_i$ for the column vector $[Z_i\edot\ep_1,\, Z_i\edot\ep_2]^t$.
Then the determinants  
$\det(\gz_i,\gz_{i+1})$ for $i=1,\ldots,p-1$ and $\det(\gz_p,\gz_1)$ are 
not negative. Moreover, $\gz_i\neq -\gz_{i+1}$ for $i=1,\ldots,p-1$ and $\gz_p\neq\pm\gz_1$.
\item
Every ordered pair $(Z_i,\,Z_j)$ of zigzags has the same orientation at all its intersection points. This is equivalent with:
$$
\forall i,j\,:\qquad |Z_i\edot Z_j|\,=\,\sharp(Z_i\cap Z_j)\,.
$$
\end{enumerate}
Condition 7 means that the homology classes are ordered counterclockwise with increasing indices and that for every homology class the indices of the zigzags in that class form a connected interval in the index set $\{1,\ldots,p\}$.
In \ref{def:very good} we formulate a condition which also restrains the ordering of zigzags within their homology class.

\

\ssnnl{Remark.}\label{rem:zigzag dimer}
The dimer model for a pattern of zigzags is the graph with a node $\bullet$ for every $+$-cell and a node $\circ$ for every $-$-cell. Two nodes are connected by an edge if the cells have a vertex in common. Figure \ref{fig:zigzag dimers}
shows an example.

\

\ssnnl{}
To a pattern of zigzags $\cZ=(Z_1,\ldots,Z_p)$ we assign the 
$2\times p$-matrix 
$$
B_\cZ\,=\,\left[
\begin{array}{rrr}
Z_1\edot\ep_1&,\ldots,&Z_p\edot\ep_1\\
Z_1\edot\ep_2&,\ldots,&Z_p\edot\ep_2
\end{array}\right]\,,
$$
displaying the intersection numbers of the zigzags in the pattern
with the two curves $\ep_1$ and $\ep_2$.
We assume that all intersections of zigzags with $\ep_1$ and $\ep_2$ are transverse.
The intersection numbers are visible in the pictures as follows.
Represent $\TT$ by the unit square with opposite sides identified.
Then $Z_j\edot\ep_2$ is the number of times the zigzag $Z_j$ crosses the right-hand vertical edge from left to right
minus the number of times it crosses from right to left.
Similarly, $Z_j\edot\ep_1$ is the number of times $Z_j$ crosses the top horizontal edge downwards minus the number of times it crosses upwards.
\\
Condition \ref{def:zigzag pattern}.1. is equivalent with
$$ 
\textit{the rank of $B_\cZ$ is $2$ and the sum of its columns is $0$.}
$$

\

\ssnnl{}
From (\ref{eq:Zee}) we see:
$$
Z_i\edot Z_j=
(Z_i\edot\ep_1)(Z_j\edot\ep_2)-(Z_i\edot\ep_2)(Z_j\edot\ep_1)=
\det\left[\!\!\begin{array}{rr} 
(Z_i\edot\ep_1)\!\!&(Z_j\edot\ep_1)\\
(Z_i\edot\ep_2)\!\!&(Z_j\edot\ep_2)
\end{array}\!\!\right].
$$
This together with Condition \ref{def:zigzag pattern}.2. implies
$$ 
Z_i\edot Z_j=0\qquad\Longleftrightarrow\qquad
[Z_i]\,=\,\pm [Z_j]\,.
$$
And thus, $\rank B_\cZ\,=\,2$ if and only if
$[Z_i]\,\neq\,\pm [Z_j]$ for some $i,j$.

\

\ssnnl{}\label{vectors cells}
Let $\cZ=\{Z_1,\ldots,Z_p\}$ be a pattern of zigzags on $\TT$. Pick a point $\star$ in one of the $-$-cells. To every $2$-cell $c$ we associate a row vector in $\ZZ^p$ as follows. Take any path $\gc$ on $\TT$ starting at $\star$
and ending in (the interior of) $c$, such that $\gc$ intersects zigzags 
transversely. Each point in $Z_j\cap\gc$ contributes, depending on the orientation, $+1$ or $-1$ to the intersection number $Z_j\edot\gc$.
Then to $c$ we associate the \textit{vector of intersection numbers}, or briefly
\textit{intersection vector}, 
$[Z_1\edot\gc,\ldots,Z_p\edot\gc]$. Choosing another path $\gc'$ from $\star$
to $c$ changes this vector by a $\ZZ$-linear combination of the rows of the matrix $B_\cZ$.

It follows from Condition \ref{zigzag pattern}.5 that at an intersection point of zigzags $Z_i$ and $Z_j$ the vectors for the two cells with unoriented boundary and the cell with negatively oriented boundary can be obtained from the vector for the $+$-cell by subtracting $1$ from the $i$-th coordinate, respectively
$1$ from the $j$-th coordinate, respectively $1$ from both $i$-th and $j$-th coordinate. We can thus capture all relevant information of the pattern of zigzags $\cZ$ in the matrix $B_\cZ$ and two additional matrices
$I_\cZ$ and $P_\cZ$, defined as follows. 

\

\ssnnl{Definition.}\label{def:IZPZQZ}
The columns of $I_\cZ$ and $P_\cZ$ correspond with the zigzags $Z_1,\ldots,Z_p$. The rows of $I_\cZ$ and $P_\cZ$ correspond with the intersection points of pairs of zigzags in $\cZ$.
The row of matrix $I_\cZ$ for a point $x\in Z_i\cap Z_j$  has $1$ in positions $i$ and $j$ and $0$ elsewhere. Matrix $P_\cZ$ has in the row for intersection point $x$ an intersection vector of the $+$-cell which has $x$ in its boundary; see Figures \ref{fig:zigzag matrices} and \ref{fig:start pattern} for examples.

It is also convenient to have the short notation $Q_\cZ=P_\cZ-I_\cZ$.
Then matrix $Q_\cZ$ has in the row for an intersection point $x$ an intersection vector of the $-$-cell which has $x$ in its boundary.

\

\ssnnl{Remark.}
Due to the $\ZZ^2 B_\cZ$-ambiguity in the choice of the intersection vectors
there is also a $\ZZ^2 B_\cZ$-ambiguity in the rows of the matrices $P_\cZ$ and $Q_\cZ$ in Definition \ref{def:IZPZQZ}.
In the algorithm we start with well-defined matrices $P_\cZ$ and $Q_\cZ$. 
In the course of the algorithm we only delete 
rows and perform the same operations on the columns of the matrices $B_\cZ$, 
$P_\cZ$, $Q_\cZ$ and $I_\cZ$ simultaneously. So the algorithm is also unambiguous. It does however happen that rows of $P_\cZ$ (resp. $Q_\cZ$)
which corresond to the same $+$-cell (resp. $-$-cell) are not the same, but differ by a vector in $\ZZ^2 B_\cZ$.

In Section \ref{sec:KZ} we pass to $\modquot{\ZZ^p}{\ZZ^2B_\cZ}$.

\

\ssnnl{Notation.}\label{notation: matrix}
For a matrix $M$ we denote $j$-th column as $M(:,j)$,
the $i$-th row as $M(i,:)$ and the $(i,j)$-entry as $M(i,j)$.

\

\ssnnl{Definition.}\label{def:pair}
Let $Z_j$ and $Z_k$ be two zigzags such that $[Z_j]\,=\,-[Z_k]$.
We say that $(Z_j,\,Z_k)$ is a \textit{$+$-opposite pair} (resp. \textit{$-$-opposite pair}) if 
$$
Q_\cZ(:,j)=-Q_\cZ(:,k)
\quad
(\textrm{resp.}\;P_\cZ(:,j)=-P_\cZ(:,k)).
$$
Suppose $Z_j\cap Z_k=\emptyset$. Then 
$(Z_j,\,Z_k)$ is a $+$-opposite pair (resp. $-$-opposite pair) if and only if
there are between
$Z_j$ and $Z_k$ no $-$-cells (resp. no $+$-cells).
Most pictures in this note contain examples of opposite pairs. 
The term `opposite pair' without $\pm$ was introduced in \cite{G} \S5.3.

\

\ssnnl{Definition.}\label{def:very good}
A good
pattern of zigzags $\cZ=(Z_1,\ldots,Z_p)$ is said to be \textit{very good} 
if it satisfies:
\begin{enumerate}\addtocounter{enumi}{8}
\item
For every homology class $[Z]$ the sequence of all zigzags
$(Z_i,\ldots,Z_{i+r}\!)$ in homology class $[Z]$ and the 
 sequence of all zigzags
$(Z_j,\ldots,Z_{j+s})$ in homology class $-[Z]$ satisfy:\\ 
\begin{equation}\label{eq:bands}
\hspace{-2.5em}
\begin{array}{ll}
\hspace{2em}\textit{for } 0\leq t\leq\min(r,s):&\hspace{-1em}
(Z_{i+t},Z_{j+t})\textit{ is a $+$-opposite pair}
\\
\hspace{2em}\textrm{and either:}&
\\
\hspace{3.5em}\textit{for } 0\leq t<\min(r,s):&\hspace{-1em}
(Z_{i+t+1},Z_{j+t})\textit{ is a $-$-opposite pair}
\\
\hspace{2em}\textrm{or: }
\textit{for } 0\leq t<\min(r,s):&\hspace{-1em}
(Z_{i+t},Z_{j+t+1})\textit{ is a $-$-opposite pair}

\end{array}
\end{equation}
\end{enumerate}
\section{The moves in the algorithm}\label{section:algorithm1}

Gulotta's algorithm transforms in an iterative way a very good pattern of zigzags
$\cZ$ into another very good one $\cZ'$.
In \cite{G} the algorithm is mainly described by transforming a drawing of  $\cZ$ into a drawing of $\cZ'$.
We will present the same algorithm by row and column operations on the matrices
$B_\cZ$, $I_\cZ$, $P_\cZ$.

\

\begin{figure}[ht]
\begin{picture}(200,50)(-50,0)
\put(50,0){\vector(0,1){40}}
\put(30,20){\vector(1,0){40}}
\put(35,28){$+$}
\put(57,8){$-$}
\put(100,20){$\rightsquigarrow$}
\put(165,0){\line(0,1){10}}
\put(165,30){\vector(0,1){10}}
\put(145,20){\line(1,0){10}}
\put(175,20){\vector(1,0){10}}
\put(155,20){\line(1,1){10}}
\put(165,10){\line(1,1){10}}
\put(150,28){$+$}
\put(172,8){$-$}
\end{picture}
\caption{\label{fig:merge}
\textit{Merging move.}}
\end{figure}
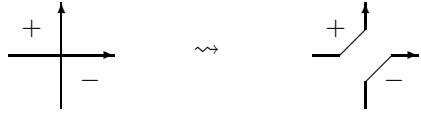

\ssnnl{Merging move.}\label{merge}
\textit{The basic move in the algorithm merges two zigzags $Z_i$ and $Z_j$
which intersect in exactly one point, as shown in 
Figure \ref{fig:merge}.} 
$\qquad$  It is evident that the merging moves preserve Conditions 1-6 in
\ref{def:zigzag pattern}.

It was pointed out in \cite{G}, that when two zigzags $Z_i$ and $Z_j$ of a pattern $\cZ$ are merged and become one zigzag $Z$ in $\cZ'$ then
\begin{equation}\label{eq:merge}
\left[\begin{array}{l} Z\edot\ep_1\\Z\edot\ep_2\end{array}\right]
\:=\:
\left[\begin{array}{l} Z_i\edot\ep_1\\Z_i\edot\ep_2\end{array}\right]
\,+\,
\left[\begin{array}{l} Z_j\edot\ep_1\\Z_j\edot\ep_2\end{array}\right]\,.
\end{equation}
Actually this means
$
[Z]=[Z_i]+[Z_j]
$
and, hence, also the column of the matrix 
$P_{\cZ'}$ which corresponds with the zigzag $Z$ is the sum of the $i$-th and the $j$-th columns of $P_\cZ$.
Moreover, as the picture indicates, the point of intersection $Z_i\cap Z_j$
disappears. The following statement also specifies where we put the new zigzag $Z$ in the list of zigzags for $\cZ'$.\\
\textbf{Conclusion:}
\textit{The merging of $Z_i$ and $Z_j$ for $Z_i\edot Z_j=1$ is given by the same column operation on $B_\cZ$, $I_\cZ$, $P_\cZ$, namely: add the $j$-th column to the 
$i$-th and subsequently delete the $j$-th column. It also deletes from $I_\cZ$ and $P_\cZ$ 
the row for $Z_i\cap Z_j$.}

\

Merging moves performed on a \textit{very good} pattern of zigzags need not preserve Conditions \ref{def:good}.7-8 and \ref{def:very good}.9. Some repairing may be needed in order to turn the pattern of zigzags produced by the merging moves into a very good one again.

\

\begin{figure}[ht]
\begin{picture}(200,50)(-50,-5)
\footnotesize

\put(10,30){\vector(1,1){10}}
\put(20,20){\line(1,1){10}}
\put(30,10){\line(1,1){10}}
\put(50,10){\vector(-1,-1){10}}
\put(30,30){\vector(1,0){20}}
\put(30,10){\vector(-1,0){20}}
\put(20,0){\line(0,1){20}}
\put(40,20){\line(0,1){20}}

\put(25,17){$-$}
\put(5,0){$+$}
\put(45,35){$+$}

\put(-2,25){$Z_i$}
\put(-2,10){$Z_j$}
\put(52,25){$Z_i$}
\put(52,10){$Z_j$}
\put(100,20){$\rightsquigarrow$}

\put(160,30){\vector(1,1){10}}
\put(200,10){\vector(-1,-1){10}}
\put(170,0){\vector(1,1){30}}
\put(190,40){\vector(-1,-1){30}}

\put(148,25){$Z_i$}
\put(148,10){$Z_j$}
\put(202,25){$Z_i$}
\put(202,10){$Z_j$}

\put(173,15){$+$}

\normalsize
\end{picture}
\caption{\label{fig:repair 1}
\textit{Repairing move 1.}}
\end{figure}
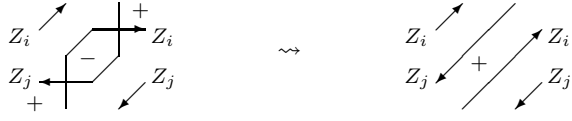

\ssnnl{Repairing move 1.}\label{repairing 1}
\textit{The first type of repairing move is shown in Figure \ref{fig:repair 1}.}
This is used when $[Z_i]\,=\,-[Z_j]$ and $\sharp(Z_i\cap Z_j)=2$,
while the area between the zigzags is just one $-$-cell (as suggested in the picture) or one $+$-cell (interchange $+$ and $-$ in the picture).
From this picture one immediately comes to the conclusion:
\\
\textbf{ Conclusion:} 
\textit{Repairing move 1 just deletes the rows for the two points of
$Z_i\cap Z_j$ from $I_\cZ$ and $P_\cZ$.}

\

\begin{figure}[ht]
\begin{picture}(200,110)(-40,-10)

\put(0,20){\vector(1,0){40}}
\put(20,0){\vector(0,1){40}}
\put(60,60){\vector(1,0){40}}
\put(80,40){\vector(0,1){40}}
\multiput(46,23)(20,10){2}{\line(2,1){10}}
\multiput(26,43)(20,10){2}{\line(2,1){10}}

\put(-12,20){$Z_i$}
\put(102,60){$Z_j$}
\put(7,0){$Z_j$}
\put(82,80){$Z_i$}

\put(5,25){$+$}
\put(25,5){$-$}
\put(65,65){$+$}
\put(85,45){$-$}

\put(130,50){$\rightsquigarrow$}

\put(160,0){
\begin{picture}(100,90)(0,0)
\put(0,20){\line(1,0){10}}
\put(30,20){\vector(1,0){10}}
\put(20,0){\line(0,1){10}}
\put(20,30){\vector(0,1){10}}
\put(60,60){\line(1,0){10}}
\put(90,60){\vector(1,0){10}}
\put(80,40){\line(0,1){10}}
\put(80,70){\vector(0,1){10}}
\put(10,20){\line(1,1){10}}
\put(20,10){\line(1,1){10}}
\put(70,60){\line(1,1){10}}
\put(80,50){\line(1,1){10}}
\multiput(46,23)(20,10){2}{\line(2,1){10}}
\multiput(26,43)(20,10){2}{\line(2,1){10}}

\put(-12,20){$Z_i$}
\put(102,60){$Z_j$}
\put(7,0){$Z_j$}
\put(82,80){$Z_i$}

\put(5,25){$+$}
\put(25,5){$-$}
\put(65,65){$+$}
\put(85,45){$-$}
\end{picture}
}
\end{picture}
\caption{\label{fig:repair 2}
\textit{Repairing move 2.}}
\end{figure}
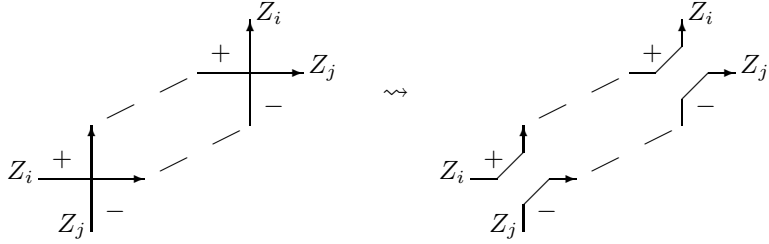

\ssnnl{Repairing move 2.}\label{repairing 2}
\textit{The second type of repairing move is shown in Figure 
\ref{fig:repair 2}.}
This is used when $[Z_i]\,=\,[Z_j]$ and $\sharp(Z_i\cap Z_j)=2$.
In this case one may distinguish three kinds of rows in the matrix $P_\cZ$,
according to whether the $i$-th entry minus the $j$-th entry is equal to $m-2$, $m-1$ or $m$, for some integer $m$ (depending on $i$ and $j$).
The rows of the latter kind correspond in the picture with 
$+$-cells in the area between the two zigzags.
When the intervals of $Z_i$ and $Z_j$ between the points of $Z_i\cap Z_j$
are swopped (as suggested by the right-hand picture) one must add $1$ to the $j$-th coordinate and subtract $1$ from the $i$-th coordinate in all rows of $P_\cZ$ corresponding with a $+$-cell in the area between $Z_i$ and $Z_j$.
One must also interchange the $i$-th and $j$-th entries
in the rows of $I_\cZ$ which correspond with intersections with the intervals of $Z_i$ and $Z_j$ between the two points of $Z_i\cap Z_j$. And one must
delete from $I_\cZ$ and $P_\cZ$ the two rows corresponding to the two points of $Z_i\cap Z_j$.
\\
\textbf{ Conclusion:}\textit{ Repairing move 2 operates on the columns of $P_\cZ$ as follows: Write $H(r)=P_\cZ(r,i)-P_\cZ(r,j)$ and $m=\max_r(H(r))$. Then
$$
\begin{array}{lll}
P_\cZ(r,i)\;\rightsquigarrow\;P_\cZ(r,i)-1\,,&
P_\cZ(r,j)\;\rightsquigarrow\;P_\cZ(r,j)+1&\textrm{if} \quad H(r)=m\,,
\\
P_\cZ(r,i)\;\rightsquigarrow\;P_\cZ(r,i)\,,&
P_\cZ(r,j)\;\rightsquigarrow\;P_\cZ(r,j)&\textrm{if} \quad H(r)\neq m\,.
\end{array}
$$
It operates on the columns of $I_\cZ$ by: 
$$
\begin{array}{lll}
I_\cZ(r,i)\;\rightsquigarrow\;I_\cZ(r,j)\,,&
I_\cZ(r,j)\;\rightsquigarrow\;I_\cZ(r,i)&
\\
&\hspace{-9em}\textrm{if} \quad H(r)=m\textrm{ and }I_\cZ(r,i)=1,
\;\textrm{or}&\hspace{-3em} H(r)=m-1 \textrm{ and } I_\cZ(r,j)=1,
\\
I_\cZ(r,i)\;\rightsquigarrow\;I_\cZ(r,i)\,,&
I_\cZ(r,j)\;\rightsquigarrow\;I_\cZ(r,j)&\textrm{otherwise}\,.
\end{array}
$$
Finally, it deletes from $I_\cZ$ and $P_\cZ$ the rows for the points of 
$Z_i\cap Z_j$.}

\begin{figure}[ht]
\begin{picture}(80,185)(-65,25)
\put(0,120){
\dashbox{1.0}(80,80){
\begin{picture}(80,80)(0,0)
\put(10,80){\vector(0,-1){80}}
\put(50,80){\vector(0,-1){80}}
\put(30,0){\vector(0,1){80}}
\put(70,0){\vector(0,1){80}}
\put(80,70){\vector(-1,0){80}}
\put(80,30){\vector(-1,0){80}}
\put(0,50){\vector(1,0){80}}
\put(0,10){\vector(1,0){80}}

\footnotesize{
\put(7,85){$1$}
\put(27,85){$2$}
\put(47,85){$3$}
\put(67,85){$4$}
\put(82,69){$1'$}
\put(82,49){$2'$}
\put(82,29){$3'$}
\put(82,9){$4'$}

\put(55,60){$+$}
\put(55,20){$+$}
\put(15,60){$+$}
\put(15,20){$+$}
\put(2,2){$-$}
\put(2,40){$-$}
\put(2,72){$-$}
\put(37,2){$-$}
\put(37,40){$-$}
\put(37,72){$-$}
\put(72,2){$-$}
\put(72,40){$-$}
\put(72,72){$-$}
}
\end{picture}
}}

\put(130,120){
\dashbox{1.0}(80,80){
\begin{picture}(80,80)(0,0)
\put(10,80){\vector(-1,-1){10}}
\put(10,60){\vector(0,-1){60}}
\put(50,80){\line(0,-1){40}}
\put(50,20){\vector(0,-1){20}}
\put(30,60){\vector(0,1){20}}
\put(30,60){\line(-1,-1){10}}
\put(30,40){\line(1,1){10}}
\put(30,0){\line(0,1){40}}
\put(70,20){\vector(0,1){60}}
\put(70,0){\vector(1,1){10}}
\put(20,70){\line(-1,-1){10}}
\put(80,70){\line(-1,0){60}}
\put(80,30){\line(-1,0){20}}
\put(40,30){\vector(-1,0){40}}
\put(0,50){\line(1,0){20}}
\put(40,50){\vector(1,0){40}}
\put(0,10){\line(1,0){60}}
\put(60,10){\line(1,1){10}}
\put(50,20){\line(1,1){10}}
\put(40,30){\line(1,1){10}}

\footnotesize{
\put(7,85){$1''$}
\put(27,85){$2''$}
\put(47,85){$3''$}
\put(67,85){$4''$}
\put(82,69){$1''$}
\put(82,49){$2''$}
\put(82,29){$3''$}
\put(82,9){$4''$}

\put(55,60){$+$}
\put(55,20){$+$}
\put(15,60){$+$}
\put(15,20){$+$}
\put(2,2){$-$}
\put(2,40){$-$}
\put(0,76){$-$}
\put(37,2){$-$}
\put(37,40){$-$}
\put(37,72){$-$}
\put(74,1){$-$}
\put(72,40){$-$}
\put(72,72){$-$}
}
\end{picture}
}}

\put(0,20){
\dashbox{1.0}(80,80){
\begin{picture}(80,80)(0,0)
\put(10,80){\vector(-1,-1){10}}
\put(10,60){\line(0,-1){20}}
\put(10,20){\vector(0,-1){20}}
\put(10,20){\line(1,1){10}}
\put(10,40){\vector(-1,-1){10}}
\put(50,60){\line(0,-1){20}}
\put(50,20){\vector(0,-1){20}}
\put(30,60){\vector(0,1){20}}
\put(30,60){\line(-1,-1){10}}
\put(30,40){\line(1,1){10}}
\put(30,20){\line(0,1){20}}
\put(70,60){\vector(0,1){20}}
\put(70,60){\line(-1,-1){10}}
\put(70,20){\line(0,1){20}}
\put(70,0){\vector(1,1){10}}
\put(20,70){\line(-1,-1){10}}
\put(80,70){\line(-1,0){20}}
\put(40,70){\line(-1,0){20}}
\put(40,70){\line(1,1){10}}
\put(60,70){\line(-1,-1){10}}
\put(80,30){\line(-1,0){20}}
\put(40,30){\line(-1,0){20}}
\put(0,50){\line(1,0){20}}
\put(40,50){\line(1,0){20}}
\put(70,40){\vector(1,1){10}}
\put(0,10){\line(1,0){20}}
\put(40,10){\line(1,0){20}}
\put(20,10){\line(1,1){10}}
\put(40,10){\line(-1,-1){10}}
\put(60,10){\line(1,1){10}}
\put(50,20){\line(1,1){10}}
\put(40,30){\line(1,1){10}}

\footnotesize{
\put(7,85){$1'''$}
\put(27,85){$2'''$}
\put(47,85){$3'''$}
\put(67,85){$4'''$}
\put(82,69){$1'''$}
\put(82,49){$2'''$}
\put(82,29){$3'''$}
\put(82,9){$4'''$}

\put(55,60){$+$}
\put(55,20){$+$}
\put(15,60){$+$}
\put(15,20){$+$}
\put(0,2){$-$}
\put(0,37){$-$}
\put(0,76){$-$}
\put(37,2){$-$}
\put(37,37){$-$}
\put(37,76){$-$}
\put(74,1){$-$}
\put(72,37){$-$}
\put(72,76){$-$}
}
\end{picture}
}}

\put(130,20){
\dashbox{1.0}(80,80){
\begin{picture}(80,80)(0,0)
\put(10,80){\vector(-1,-1){10}}
\put(50,80){\vector(-1,-1){50}}
\put(80,70){\vector(-1,-1){70}}
\put(80,30){\vector(-1,-1){30}}
\put(30,0){\vector(1,1){50}}
\put(70,0){\vector(1,1){10}}
\put(0,50){\vector(1,1){30}}
\put(0,10){\vector(1,1){70}}

\footnotesize{
\put(7,85){$1'''''$}
\put(27,85){$2'''$}
\put(47,85){$3''''$}
\put(67,85){$4''''$}
\put(82,69){$1''''$}
\put(82,49){$2''''$}
\put(82,29){$3''''$}
\put(82,9){$4''''$}

\put(60,20){$-$}
\put(40,40){$-$}
\put(20,60){$-$}
\put(74,1){$-$}
\put(0,76){$-$}
}
\end{picture}
}}

\end{picture}
\caption{\label{fig:band}
\textit{Transforming an alternating sequence of opposite pairs.
}}
\end{figure}
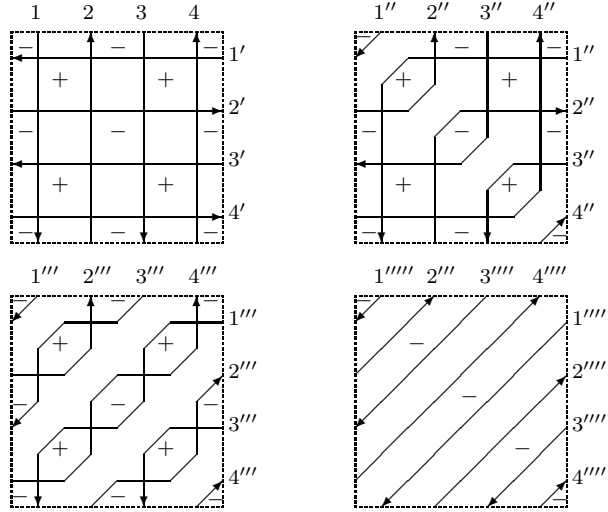

\

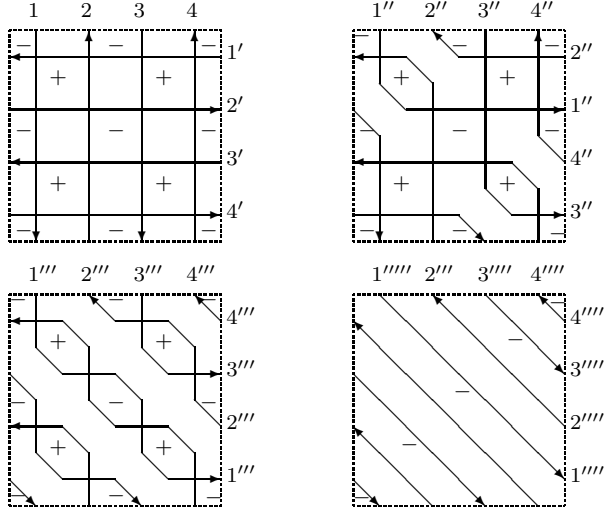
\begin{figure}[t]
\begin{picture}(80,185)(-65,25)
\put(0,120){
\dashbox{1.0}(80,80){
\begin{picture}(80,80)(0,0)
\put(10,80){\vector(0,-1){80}}
\put(50,80){\vector(0,-1){80}}
\put(30,0){\vector(0,1){80}}
\put(70,0){\vector(0,1){80}}
\put(80,70){\vector(-1,0){80}}
\put(80,30){\vector(-1,0){80}}
\put(0,50){\vector(1,0){80}}
\put(0,10){\vector(1,0){80}}

\footnotesize{
\put(7,85){$1$}
\put(27,85){$2$}
\put(47,85){$3$}
\put(67,85){$4$}
\put(82,69){$1'$}
\put(82,49){$2'$}
\put(82,29){$3'$}
\put(82,9){$4'$}

\put(55,60){$+$}
\put(55,20){$+$}
\put(15,60){$+$}
\put(15,20){$+$}
\put(2,2){$-$}
\put(2,40){$-$}
\put(2,72){$-$}
\put(37,2){$-$}
\put(37,40){$-$}
\put(37,72){$-$}
\put(72,2){$-$}
\put(72,40){$-$}
\put(72,72){$-$}
}
\end{picture}
}}

\put(130,120){
\dashbox{1.0}(80,80){
\begin{picture}(80,80)(0,0)
\put(10,80){\line(0,-1){20}}
\put(10,40){\vector(0,-1){40}}
\put(50,80){\line(0,-1){60}}
\put(30,0){\line(0,1){60}}
\put(70,0){\line(0,1){20}}
\put(70,40){\vector(0,1){40}}
\put(20,70){\vector(-1,0){20}}
\put(80,70){\line(-1,0){40}}
\put(60,30){\vector(-1,0){60}}
\put(20,50){\vector(1,0){60}}
\put(0,10){\line(1,0){40}}
\put(60,10){\vector(1,0){20}}
\put(40,70){\vector(-1,1){10}}
\put(30,60){\line(-1,1){10}}
\put(20,50){\line(-1,1){10}}
\put(10,40){\line(-1,1){10}}
\put(40,10){\vector(1,-1){10}}
\put(50,20){\line(1,-1){10}}
\put(60,30){\line(1,-1){10}}
\put(70,40){\line(1,-1){10}}

\footnotesize{
\put(7,85){$1''$}
\put(27,85){$2''$}
\put(47,85){$3''$}
\put(67,85){$4''$}
\put(82,69){$2''$}
\put(82,49){$1''$}
\put(82,29){$4''$}
\put(82,9){$3''$}

\put(55,60){$+$}
\put(55,20){$+$}
\put(15,60){$+$}
\put(15,20){$+$}
\put(2,2){$-$}
\put(2,40){$-$}
\put(0,76){$-$}
\put(37,2){$-$}
\put(37,40){$-$}
\put(37,72){$-$}
\put(74,1){$-$}
\put(72,40){$-$}
\put(72,72){$-$}
}
\end{picture}
}}

\put(0,20){
\dashbox{1.0}(80,80){
\begin{picture}(80,80)(0,0)
\put(10,80){\line(0,-1){20}}
\put(10,40){\line(0,-1){20}}
\put(50,80){\line(0,-1){20}}
\put(50,40){\line(0,-1){20}}
\put(30,0){\line(0,1){20}}
\put(30,40){\line(0,1){20}}
\put(70,0){\line(0,1){20}}
\put(70,40){\line(0,1){20}}
\put(20,70){\vector(-1,0){20}}
\put(60,70){\line(-1,0){20}}
\put(60,30){\line(-1,0){20}}
\put(20,30){\vector(-1,0){20}}
\put(20,50){\line(1,0){20}}
\put(60,50){\vector(1,0){20}}
\put(20,10){\line(1,0){20}}
\put(60,10){\vector(1,0){20}}
\put(40,70){\vector(-1,1){10}}
\put(30,60){\line(-1,1){10}}
\put(20,50){\line(-1,1){10}}
\put(10,40){\line(-1,1){10}}
\put(40,10){\vector(1,-1){10}}
\put(50,20){\line(1,-1){10}}
\put(60,30){\line(1,-1){10}}
\put(70,40){\line(1,-1){10}}
\put(0,10){\vector(1,-1){10}}
\put(10,20){\line(1,-1){10}}
\put(20,30){\line(1,-1){10}}
\put(30,40){\line(1,-1){10}}
\put(40,50){\line(1,-1){10}}
\put(50,60){\line(1,-1){10}}
\put(60,70){\line(1,-1){10}}
\put(80,70){\vector(-1,1){10}}

\footnotesize{
\put(7,85){$1'''$}
\put(27,85){$2'''$}
\put(47,85){$3'''$}
\put(67,85){$4'''$}
\put(82,69){$4'''$}
\put(82,49){$3'''$}
\put(82,29){$2'''$}
\put(82,9){$1'''$}

\put(55,60){$+$}
\put(55,20){$+$}
\put(15,60){$+$}
\put(15,20){$+$}
\put(0,2){$-$}
\put(0,37){$-$}
\put(0,76){$-$}
\put(37,2){$-$}
\put(37,37){$-$}
\put(37,76){$-$}
\put(74,1){$-$}
\put(72,37){$-$}
\put(72,76){$-$}
}
\end{picture}
}}

\put(130,20){
\dashbox{1.0}(80,80){
\begin{picture}(80,80)(0,0)
\put(10,80){\vector(1,-1){70}}
\put(50,80){\vector(1,-1){30}}
\put(0,50){\vector(1,-1){50}}
\put(0,10){\vector(1,-1){10}}
\put(30,0){\vector(-1,1){30}}
\put(70,0){\vector(-1,1){70}}
\put(80,30){\vector(-1,1){50}}
\put(80,70){\vector(-1,1){10}}

\footnotesize{
\put(7,85){$1'''''$}
\put(27,85){$2'''$}
\put(47,85){$3''''$}
\put(67,85){$4''''$}
\put(82,69){$4''''$}
\put(82,49){$3''''$}
\put(82,29){$2''''$}
\put(82,9){$1''''$}

\put(0,1){$-$}
\put(18,21){$-$}
\put(38,41){$-$}
\put(58,61){$-$}
\put(74,75){$-$}
}
\end{picture}
}}

\end{picture}
\caption{\label{fig:band2}
\textit{Transforming an alternating sequence of opposite pairs.
}}
\end{figure}

\ssnnl{Example.}\label{example:band}
Let $(Z_1,\ldots,Z_{2s})$ and $(Z'_1,\ldots,Z'_{2s})$
be two sequences of zigzags in a pattern $\cZ^{(1)}$ such that 
$\sharp(Z_i\cap Z'_j)=1$ for all $i, j$ and $Z_i\cap Z_j\,=\,Z'_i\cap Z'_j=\emptyset$ for all $i\neq j$
and such that $(Z_{j-1},Z_j)$
and $(Z'_{j-1},Z'_j)$ are $(-1)^j$-opposite pairs
for $j=2,\ldots,2s$. There are two well controllable cases in which merging of these two sequences followed by repairing moves 2 and 1 yields a sequence of
zigzags $(Z''''_1,\ldots,Z''''_{2s})$ such that
$Z''''_i\cap Z''''_j=\emptyset$ for all $i\neq j$ and such that
$(Z''''_{j-1},Z''''_j)$ is a $(-1)^j$-opposite pair
for $j=2,\ldots,2s$. 
\begin{eqnarray*}
\textbf{Case 1:}&&
\textit{$Z_j$ and $Z'_j$ merge for $j=1,\ldots,2s$.}
\\
\textbf{Case 2:}&&
\textit{$Z_j$ and $Z'_{j-(-1)^j}$ merge for $j=1,\ldots,2s$.}
\end{eqnarray*}
Figures \ref{fig:band} and \ref{fig:band2} show this for $s=2$ and clearly generalize to arbitrary $s$.

In either case let $(Z''_1,\ldots,Z''_{2s})$ be the sequence of zigzags in the pattern $\cZ^{(2)}$  
which results from the merging.
Now transform $\cZ^{(2)}$ by applying repairing moves of type 2 
at the points of $Z''_i\cap Z''_j$ for all $i\neq j$ for which $[Z''_i]=[Z''_j]$.
The result is the sequence of zigzags  $(Z'''_1,\ldots,Z'''_{2s})$ in the
pattern $\cZ^{(3)}$.
Then $Z'''_i\cap Z'''_j\neq \emptyset$ only if $[Z'''_i]=-[Z'''_j]$.
Next transform $\cZ^{(3)}$ by applying repairing moves of type 1 
at the points of $Z'''_i\cap Z'''_j$ for all $i,j$.
Call the resulting pattern of zigzags $\cZ^{(4)}$ and the relevant sequence of zigzags $(Z''''_1,\ldots,Z''''_{2s})$.
In this last sequence $Z''''_i\cap Z''''_j=\emptyset$ 
for all $i\neq j$ and $(Z''''_{j-1},Z''''_j)$
is a $(-1)^j$-opposite pair for $j=2,\ldots,2s$.

It is instructive to perform the moves in Figures \ref{fig:band} and
\ref{fig:band2} also for the
matrices $B_\cZ$, $I_\cZ$, $P_\cZ$. For the top-left picture 
$B_\cZ$, $I_\cZ$, $P_\cZ$ are given in Figure \ref{fig:start pattern}.
For the other pictures one may follow the description of the merging 
and repairing moves.

\

\ssnnl{Remark.}\label{remark:reorder}
In \ref{example:band} we have chosen the labels for the zigzags while drawing the pictures. If one uses the matrix operations instead, the algorithm determines the labels and it may be necessary to reorder the columns of the matrices
$P_\cZ$, $Q_\cZ$, $I_\cZ$ to meet the requirements of Equation 
(\ref{eq:bands}). Example \ref{example:band} shows that such a reordering is always possible.

\

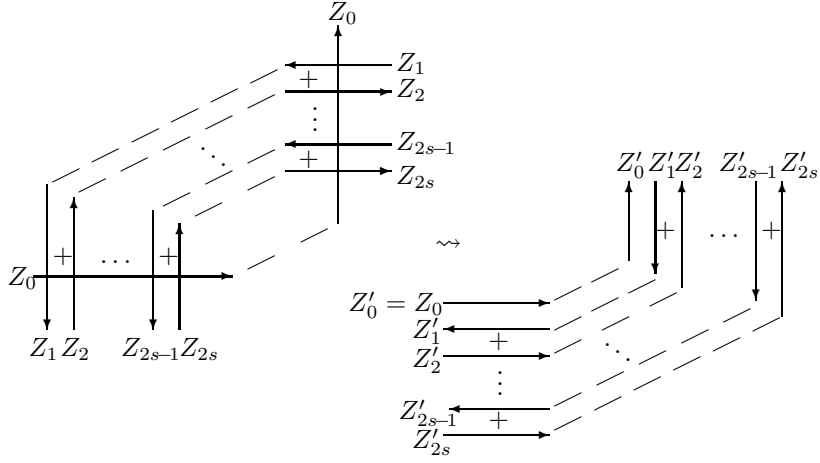
\begin{figure}[t]
\begin{picture}(200,175)(-60,-10)
\put(5,30){
\begin{picture}(100,90)(0,0)
\put(-35,20){\vector(1,0){75}}
\put(20,0){\vector(0,1){40}}
\put(60,60){\vector(1,0){40}}
\put(-20,0){\vector(0,1){50}}
\put(60,60){\vector(1,0){40}}
\put(80,40){\vector(0,1){75}}
\put(100,70){\vector(-1,0){40}}
\put(10,45){\vector(0,-1){45}}
\put(60,90){\vector(1,0){40}}
\put(100,100){\vector(-1,0){40}}
\put(-30,55){\vector(0,-1){55}}

\multiput(46,23)(20,10){2}{\line(2,1){10}}
\multiput(26,43)(20,10){2}{\line(2,1){10}}
\multiput(14,47)(16,8){3}{\line(2,1){10}}
\multiput(-18,52)(16,8){5}{\line(2,1){10}}
\multiput(-28,56)(15,7.5){6}{\line(2,1){10}}

\put(-10,25){$\ldots$}
\put(27,61){$\ddots$}
\put(70,75){$\vdots$}

\put(-45,17){$Z_0$}
\put(76,117){$Z_0$}
\put(102,55){$Z_{2s}$}
\put(102,68){$Z_{2s\!-\!1}$}
\put(20,-10){$Z_{2s}$}
\put(-2,-10){$Z_{2s\!-\!1}$}
\put(-37,-10){$Z_1$}
\put(-25,-10){$Z_2$}
\put(102,98){$Z_1$}
\put(102,88){$Z_2$}

\put(12,25){$+$}
\put(-28,25){$+$}
\put(65,62){$+$}
\put(65,92){$+$}

\end{picture}
}

\put(125,60){$\rightsquigarrow$}

\put(130,0){
\begin{picture}(100,90)(0,0)
\put(-5,20){\vector(1,0){40}}
\put(35,30){\vector(-1,0){40}}
\put(-5,40){\vector(1,0){40}}
\put(85,46){\vector(0,1){40}}
\put(75,86){\vector(0,-1){35}}
\put(65,56){\vector(0,1){30}}
\put(35,0){\vector(-1,0){38}}
\put(-5,-10){\vector(1,0){40}}
\put(123,35){\vector(0,1){51}}
\put(113,86){\vector(0,-1){45}}

\multiput(37,21)(18,9){3}{\line(2,1){10}}
\multiput(37,30)(14,7){3}{\line(2,1){10}}
\multiput(37,41)(14,7){2}{\line(2,1){10}}
\multiput(37,1)(16,8){5}{\line(2,1){10}}
\multiput(37,-8)(15,7.5){6}{\line(2,1){10}}

\put(15,5){$\vdots$}
\put(55,17){$\ddots$}
\put(95,65){$\ldots$}

\put(-17,16){$Z'_2$}
\put(-17,27){$Z'_1$}
\put(-41,37){$Z'_0=Z_0$}
\put(82,90){$Z'_2$}
\put(72,90){$Z'_1$}
\put(60,90){$Z'_0$}
\put(-17,-15){$Z'_{2s}$}
\put(-22,-4){$Z'_{2s\!-\!1}$}
\put(122,90){$Z'_{2s}$}
\put(99,90){$Z'_{2s\!-\!1}$}

\put(12,23){$+$}
\put(75,65){$+$}
\put(12,-7){$+$}
\put(115,65){$+$}

\end{picture}
}
\end{picture}

\caption{\label{fig:repair 3}
\textit{Repairing move 3.}}
\end{figure}

\ssnnl{Repairing move 3.}\label{repairing 3}
\textit{The third type of repairing move is shown in Figure \ref{fig:repair 3}.}
It is used for a zigzag $Z_0$ and a sequence of zigzags $(Z_1,\ldots,Z_{2s})$
which satisfy the following conditions.
Firstly,  
$Z_i\cap Z_j=\emptyset$ for all $i>j\geq 1$ and $(Z_{j-1},Z_j)$
is a $(-1)^j$-opposite pair for $j=2,\ldots,2s$.
Secondly $[Z_0]=\pm[Z_1]$ and 
$\sharp(Z_0\cap Z_j)=2$ for $j=1,\ldots,2s$.
Reversing if necessary the labeling in the sequence $(Z_1,\ldots,Z_{2s})$
we may without loss of generality assume $[Z_0]=(-1)^j[Z_j]$ for $j=1,\ldots,2s$.

Gulotta's instructions (cf. \cite{G} \S 5.3) in this situation are to remove
$(Z_1,\ldots,Z_{2s})$ and to insert a sequence of zigzags 
$(Z'_1,\ldots,Z'_{2s})$ such that $Z'_i\cap Z'_j=\emptyset$ 
for all $j>i\geq 0$ and such that
$(Z'_{j-1},Z'_j)$ is a $(-1)^j$-opposite pair for $j=1,\ldots,2s$. For notational convenience we write here and below $Z'_0=Z_0$.

In terms of the matrices $I_\cZ$ and $P_\cZ$ this means that
we first delete from $I_\cZ$ and $P_\cZ$ all rows which correspond with an intersection point on one of the zigzags $Z_1,\ldots,Z_{2s}$ and subsequently
replace, for $j=1,\ldots,2s$, the column of $P_\cZ$ which corresponds with the zigzag $Z_j$ by $(-1)^j$ times the column of $P_\cZ$ which corresponds with the zigzag $Z_0$.

\begin{figure}[ht]
\begin{picture}(120,55)(-100,0)
\put(0,20){\vector(1,0){120}}
\put(20,0){\vector(0,1){40}}
\put(50,0){\vector(0,1){40}}
\put(100,0){\vector(0,1){40}}
\put(35,40){\vector(0,-1){40}}
\put(85,40){\vector(0,-1){40}}

\put(7,30){$+$}
\put(37,30){$+$}
\put(87,30){$+$}
\put(22,10){$-$}
\put(60,10){$\ldots$}

\put(15,45){$Z'_0$}
\put(30,45){$Z'_1$}
\put(45,45){$Z'_2$}
\put(73,45){$Z'_{2s\!-\!1}$}
\put(97,45){$Z'_{2s}$}
\put(-15,18){$Z_\infty$}
\end{picture}
\caption{\label{fig:band3}
\textit{Intersections with alternating sequence of opposite pairs.
}}
\end{figure}
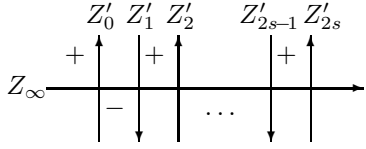

Next we expand every row of $I_\cZ$ and $P_\cZ$ which corresponds with an intersection point of $Z_0$ and a zigzag 
$Z_\infty\neq Z_0,Z_1,\ldots,Z_{2s}$ to $1+2s$ rows which correspond with the intersection points of $Z_\infty$ with $Z'_0,Z'_1,\ldots,Z'_{2s}$
(see Figure \ref{fig:band3}). The columns of $I_\cZ$ and $P_\cZ$ are labeled in such a way that the column which originally corresponded to $Z_j$ now corresponds to $Z'_j$ for $j=0,\ldots,2s$. 

Thus a row of $I_\cZ$ for an intersection point $x$ of
$Z_\infty$ and $Z_0$ is replaced by $1+2s$ rows of which the $j$-th one
(for $j=0,\ldots,2s$) has entry $1$ in the columns corresponding with the zigzags $Z_\infty$ and $Z'_j$.
Of the $1+2s$ new rows of $P_\cZ$ the $0$-th one is equal to the row
of $P_\cZ$ which corresponds to $x$. Figure \ref{fig:band3} shows that for odd $j$ the $j$-th row is obtained from the $(j-1)$-st one by adding 
$-1$ in the column for $Z'_{j-1}$ and $+1$ in the column for $Z'_j$;
for even $j\geq 2$ the $j$-th row is equal to the $(j-1)$-st one.

\
\section{Running the algorithm}\label{section:algorithm2}
In this Section we translate the algorithm described by Gulotta in terms of pictures,
into an iterative proces operating on matrices.

\

\ssnnl{}\label{BA}
In order to prepare the input for the algorithm from the set 
$\cA=\{\va_1,\ldots,\va_N\}$ (see \S \ref{intro}) we take a $2\times N$-matrix
$B_\cA$ such that its rows are a $\ZZ$-basis for the lattice
$
\{(\ell_1,\ldots,\ell_N)\in\ZZ^N\:|\:\ell_1\va_1+\ldots+\ell_N\va_N\,=\,0\}
$.

\textit{The aim of the algorithm is to create a very good pattern of zigzags $\cZ$ such that $B_\cZ$ results by permuting
and splitting up the columns of $B_\cA$ as follows}

\begin{equation}\label{eq:primitive}
\textrm{column
$\left[\!\!\begin{array}{r} r\\ s\end{array}\!\!\right]$
of $B_\cA$}
\;\rightsquigarrow\;
\textrm{$d=\mathrm{g.c.d.}(r,s)$ columns 
$\frac{1}{d}\left[\!\!\begin{array}{r} r\\ s\end{array}\!\!\right]$
in $B_\cZ$.}
\end{equation}

\

\ssnnl{}\label{start algorithm}
The algorithm starts
with a very good pattern of zigzags $\cZ$ for which $B_\cZ$ is the $2\times (2n_1+2n_2)$-matrix
$$
\raisebox{-1.4ex}{$\left[\rule{0mm}{3ex}\right.$}
\begin{array}{llll}
\overbrace{1\dots1}^{n_1}&\overbrace{0\dots0}^{n_2}&
\overbrace{-1\dots-1}^{n_1}&\overbrace{\hspace{.8em}0\dots\hspace{.9em}0}^{n_2}\\
0\dots0&1\dots1&\hspace{.8em}0\dots\hspace{.9em}0&-1\dots-1
\end{array}
\raisebox{-1.4ex}{$\left.\rule{0mm}{3ex}\right]$}\,,
$$
where $n_1$ (resp. $n_2$) is the sum of the positive entries in the first (resp. second) row of $B_\cA$.
This matrix is realized by a pattern with straight lines, $n_1$ vertically down,
$n_1$ vertically up, $n_2$ horizontally left-to-right and $n_2$ horizontally right-to-left. To get a very good pattern one takes the vertical lines alternatingly down and up, and the horizontal lines alternatingly left-to-right and right-to-left. 
The matrices $B_\cZ,\,P_\cZ$ and $I_\cZ$ for the initial pattern
have $2n_1+2n_2$ columns.
There are $4n_1n_2$ intersection points and, hence, $I_\cZ$ and $P_\cZ$
have $4n_1n_2$ rows. 
We build $I_\cZ$ and $P_\cZ$ as follows:
the $+$-cells are given by pairs $(a,b)$ in 
$\{1,\ldots,n_1\}\times\{1,\ldots,n_2\}$. 
The rows of $I_\cZ$ which correspond to the four vertices of the $+$-cell $(a,b)$ are have $1$ in positions $a,\,b+n_1$, resp.
$a,\,b+2n_1+n_2$, resp.  $a+n_1+n_2,\,b+n_1$, resp. $a+n_1+n_2,\,b+2n_1+n_2$. All other entries in these rows of $I_\cZ$ are $0$.
The non-zero entries of the four rows of $P_\cZ$ which correspond to the vertices of the $+$-cell $(a,b)$ 
are $1$ in position $j$ if $1\leq j\leq a$ or $2n_1+n_2+1\leq j\leq 2n_1+n_2+b$
and $-1$ in position $j$ if
$n_1+n_2+1\leq j \leq n_1+n_2+a-1$ or $n_1+1\leq j\leq n_1+b-1$.
Figure \ref{fig:start pattern} shows an example
with $n_1=n_2=2$. 
 
\begin{figure}[ht]
\begin{picture}(80,165)(-40,0)
\put(0,40){
\dashbox{1.0}(80,80){
\begin{picture}(80,80)(0,0)
\put(10,80){\vector(0,-1){80}}
\put(50,80){\vector(0,-1){80}}
\put(30,0){\vector(0,1){80}}
\put(70,0){\vector(0,1){80}}
\put(80,70){\vector(-1,0){80}}
\put(80,30){\vector(-1,0){80}}
\put(0,50){\vector(1,0){80}}
\put(0,10){\vector(1,0){80}}

\footnotesize{
\put(7,85){$1$}
\put(47,85){$2$}
\put(27,-7){$5$}
\put(67,-7){$6$}
\put(83,69){$7$}
\put(83,29){$8$}
\put(-10,49){$3$}
\put(-10,9){$4$}

\put(55,60){$+$}
\put(55,20){$+$}
\put(15,60){$+$}
\put(15,20){$+$}

}
\end{picture}
}}
\put(95,80){\footnotesize{
$
\begin{array}{rl}
\begin{array}{ll}
&\hspace{-.5em}B_\cZ\\
I_\cZ&\hspace{-.5em}P_\cZ
\end{array}&\hspace{-1.5em}
\left[\!\!\begin{array}{rrrrrrrr}
1\!\!&1\!\!&0\!\!&0\!\!&-1\!\!&-1\!\!&0\!\!&0\\
0\!\!&0\!\!&1\!\!&1\!\!&0\!\!&0\!\!&-1\!\!&-1
\end{array}\!\!\right]\\ [2ex]
\left[\!\!\begin{array}{l}
10100000\\
10000010\\
00101000\\
00001010\\
10010000\\
10000001\\
00011000\\
00001001\\
01100000\\
01000010\\
00100100\\
00000110\\
01010000\\
01000001\\
00010100\\
00000101
\end{array}\!\!\right]&\hspace{-1.5em}
\left[\!\!\begin{array}{rrrrrrrr}
1\!\!&0\!\!&0\!\!&0\!\!&0\!\!&0\!\!&1\!\!&0\\
1\!\!&0\!\!&0\!\!&0\!\!&0\!\!&0\!\!&1\!\!&0\\
1\!\!&0\!\!&0\!\!&0\!\!&0\!\!&0\!\!&1\!\!&0\\
1\!\!&0\!\!&0\!\!&0\!\!&0\!\!&0\!\!&1\!\!&0\\
1\!\!&0\!\!&-1\!\!&0\!\!&0\!\!&0\!\!&1\!\!&1\\
1\!\!&0\!\!&-1\!\!&0\!\!&0\!\!&0\!\!&1\!\!&1\\
1\!\!&0\!\!&-1\!\!&0\!\!&0\!\!&0\!\!&1\!\!&1\\
1\!\!&0\!\!&-1\!\!&0\!\!&0\!\!&0\!\!&1\!\!&1\\
1\!\!&1\!\!&0\!\!&0\!\!&-1\!\!&0\!\!&1\!\!&0\\
1\!\!&1\!\!&0\!\!&0\!\!&-1\!\!&0\!\!&1\!\!&0\\
1\!\!&1\!\!&0\!\!&0\!\!&-1\!\!&0\!\!&1\!\!&0\\
1\!\!&1\!\!&0\!\!&0\!\!&-1\!\!&0\!\!&1\!\!&0\\
1\!\!&1\!\!&-1\!\!&0\!\!&-1\!\!&0\!\!&1\!\!&1\\
1\!\!&1\!\!&-1\!\!&0\!\!&-1\!\!&0\!\!&1\!\!&1\\
1\!\!&1\!\!&-1\!\!&0\!\!&-1\!\!&0\!\!&1\!\!&1\\
1\!\!&1\!\!&-1\!\!&0\!\!&-1\!\!&0\!\!&1\!\!&1
\end{array}\!\!\right]
\end{array}
$
}}

\end{picture}
\caption{\label{fig:start pattern}
\textit{Pattern to start from.
}}
\end{figure}

\

\ssnnl{}\label{iteration step}
At the beginning of an iteration step we have a very good pattern of zigzags 
$\cZ=\{Z_1,\ldots,Z_p\}$.
The merging moves are determined from the positions of the columns of the matrix $B_\cA$ with respect to the columns of the 
matrix $B_\cZ$. The columns of $B_\cZ$ are vectors in the plane $\RR^2$
and the ordering by increasing index coincides with the counter-clockwise
cyclic ordering. In agreement with this cyclic structure we treat the first column of $B_\cZ$ as consecutive to the last column.

A column $v$ of $B_\cA$ is either a positive integer multiple of a column of $B_\cZ$ or there is a unique pair of consecutive columns $v_1$ and $v_2$ of $B_\cZ$ such that $v=c_1v_1+c_2v_2$ with $c_1,\,c_2\in\QQ_{>0}$.

\

\ssnnl{} The \textbf{algorithm terminates} automatically when all columns of $B_\cA$ are multiples of columns of $B_\cZ$. Since the merging moves decrease the number of zigzags the algorithm will surely terminate.

\

\ssnnl{}\label{cramer} Cramer's rule explicates the relation $v=c_1v_1+c_2v_2$:
\begin{equation}\label{eq:cramer}
\det(v_1,v_2)\,v\:=\:\det(v,v_2)\,v_1\,-\,\det(v,v_1)\,v_2\,.
\end{equation}
In the pattern we start with the determinants of consecutive pairs of 
non-equal columns of $B_\cZ$ are $1$. The merging of two zigzags in a very good pattern $\cZ$ with exactly one intersection point replaces the corresponding columns of $B_\cZ$ by their sum. Thus in the next iteration step in the algorithm Equation (\ref{eq:cramer}) becomes
$$ 
\det(v_1,v_2)\,v\:=\:
(\det(v,v_2)-m)\,v_1\,+\,m(v_1+v_2)\,-\,(\det(v,v_1)+m)\,v_2
$$
with $m\,=\,\min(\det(v_1,v),\det(v,v_2))$. This then gives $v$ either as a multiple of $v_1+v_2$ or as a positive linear combination of
$v_1$ and $v_1+v_2$ or of $v_1+v_2$ and $v_2$. 
Note that
$
\det(v_1,v_1+v_2)=\det(v_1+v_2,v_2)=\det(v_1,v_2)
$.
\\
\textbf{Conclusion:}
\textit{In all cases in which Equation (\ref{eq:cramer}) is used in the algorithm $\det(v_1,v_2)\,=\,1$ and the equation actually reads
\begin{equation}\label{eq:cramer2}
v\:=\:\det(v,v_2)\,v_1\,-\,\det(v,v_1)\,v_2\,.
\end{equation}
Column $v$ of $B_\cA$ thus leads to the command that
$m$ zigzags of the pattern $\cZ$ in the homology class corresponding with the column $v_1$ of $B_\cZ$ must merge with
$m$ zigzags in the homology class corresponding with the column $v_2$;
here $m\,=\,\min(\det(v_1,v),\det(v,v_2))$.}

\

\ssnnl{}\label{actual iteration}
As Equation (\ref{eq:cramer2}) indicates
we need the determinants of the $2\times 2$-matrices with first column from $B_\cA$
and second column from $B_\cZ$. These are simultaneously given as the entries of the $N\times p$-matrix
$$
S=B_\cA^t\,J\,B_\cZ\qquad\textrm{with}\quad
J=\left[\begin{array}{rr} 0&1\\-1&0\end{array}\right]\,.
$$
The columns of $B_\cA$ correspond with the rows of $S$. One can implement
the discussion in \ref{iteration step} and \ref{cramer} for all
columns of $B_\cA$ simultaneously as follows.
Let $S^c=S(:,[2:p,1])$
be the $N\times p$-matrix obtained from $S$ by cyclically permuting the columns so that the first column comes in the last position.
Let $R$ be the $N\times p$-matrix with $(i,j)$-entry
$$
\begin{array}{rcll}
R_{ij}&=&\textstyle{\frac{1}{2}}(|S_{ij}|+|S^c_{ij}|-|S_{ij}+S^c_{ij}|)
&\textrm{if}\;S_{ij}<0\,,
\\
R_{ij}&=&0&\textrm{if}\;S_{ij}\geq 0\,.
\end{array}
$$
We define functions $\rho$ and $\gl$ on $\{1,\ldots,p\}$ by
$$
\rho(j)\,=\,\sum_{i=1}^N R_{ij}\,,
\qquad
\gl(j)\,=\,\rho(j-1)\;\textrm{if}\;j>1\,,\quad \gl(1)\,=\,\rho(p)\,.
$$
Next we define for a homology class of zigzags $[Z]$ in the pattern
$\cZ=(Z_1,\ldots,Z_p)$:
\begin{eqnarray*}
\widetilde{\rho}([Z])&=&\max\{\rho(j)\:|\:Z_j\in[Z]\}\,,
\\
\widetilde{\gl}([Z])&=&\max\{\gl(j)\:|\:Z_j\in[Z]\}\,,
\\
\widetilde{\mu}([Z])&=&\sharp ([Z])\,-\,\widetilde{\rho}([Z])\,-\,\widetilde{\gl}([Z])\,.
\end{eqnarray*}
Then $\widetilde{\rho}([Z])$ (resp. $\widetilde{\gl}([Z])\:$) is the number of zigzags in $[Z]$
which must merge with a zigzag in the homology class immediately after
(resp. before) $[Z]$ and $\widetilde{\mu}([Z])$ is the number of zigzags in $[Z]$ which must not merge with another zigzag.
The merging step defined in \ref{merging} decreases the number of zigzags by
$$
q=\sum_{h=1}^p\gl(h)
$$
and uses the map
$\gf:\{1,\ldots,p\}\longrightarrow\{1,\ldots,p-q\}$,
\begin{equation}\label{eq:pre-merge}
\gf(j)\,=\,\left\{
\begin{array}{lll}
j-\sum_{h=1}^j\gl(h)&\textrm{if}&j>\gl(1)\\
p-q+j-\sum_{h=1}^j\gl(h)&\textrm{if}&j\leq\gl(1)\,.
\end{array}\right.
\end{equation}

\

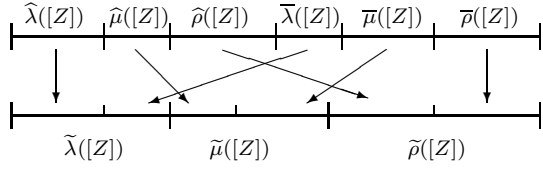
\begin{figure}[ht]
\begin{picture}(300,60)(-75,0)
\put(0,15){\line(1,0){200}}
\put(0,45){\line(1,0){200}}
\put(0,40){\line(0,1){10}}
\put(35,40){\line(0,1){10}}
\put(60,40){\line(0,1){10}}
\put(100,40){\line(0,1){10}}
\put(125,40){\line(0,1){10}}
\put(160,40){\line(0,1){10}}
\put(200,40){\line(0,1){10}}

\put(0,10){\line(0,1){10}}
\put(60,10){\line(0,1){10}}
\put(120,10){\line(0,1){10}}
\put(200,10){\line(0,1){10}}
\put(35,15){\line(0,1){4}}
\put(85,15){\line(0,1){4}}
\put(160,15){\line(0,1){4}}

\put(17,40){\vector(0,-1){20}}
\put(47,40){\vector(1,-1){20}}
\put(80,40){\vector(3,-1){55}}
\put(112,40){\vector(-3,-1){60}}
\put(142,40){\vector(-3,-2){30}}
\put(180,40){\vector(0,-1){20}}

\footnotesize
\put(5,50){$\widehat{\gl}([Z])$}
\put(37,50){$\widehat{\mu}([Z])$}
\put(68,50){$\widehat{\rho}([Z])$}
\put(102,50){$\overline{\gl}([Z])$}
\put(133,50){$\overline{\mu}([Z])$}
\put(170,50){$\overline{\rho}([Z])$}
\put(20,0){$\widetilde{\gl}([Z])$}
\put(75,0){$\widetilde{\mu}([Z])$}
\put(150,0){$\widetilde{\rho}([Z])$}
\normalsize
\end{picture}
\caption{\label{fig:shuffle}
\textit{Shuffle within one homology class}}
\end{figure}

\ssnnl{} In order to eventually satisfy Requirement \ref{def:very good}.9
and to benefit from Example \ref{example:band} we permute the zigzags in each homology class as follows.
First we define for each homology class of zigzags $[Z]$ for which the opposite class $-[Z]$ also occurs in the pattern $\cZ$:
$$
\begin{array}{rcl}
\widehat{\rho}([Z])\,=\,\min\{\widetilde{\rho}([Z]),\widetilde{\rho}(-[Z]))
\,,&\;&
\widehat{\gl}([Z])\,=\,\min\{\widetilde{\gl}([Z]),\widetilde{\gl}(-[Z]))
\,,
\\
\widehat{\mu}([Z])\,=\,\min\{\widetilde{\mu}([Z]),\widetilde{\mu}(-[Z]))
\,.&&
\end{array}
$$
If $-[Z]$ does not occur in $\cZ$ we put
$\widehat{\rho}([Z])=\widehat{\gl}([Z])=\widehat{\mu}([Z])=0$.
\\
Next we write for every homology class $[Z]$ in $\cZ$:
$$
\overline{\rho}([Z])\!=\!\widetilde{\rho}([Z])-\widehat{\rho}([Z]),\:
\overline{\gl}([Z])\!=\!\widetilde{\gl}([Z])-\widehat{\gl}([Z]),\:
\overline{\mu}([Z])\!=\!\widetilde{\mu}([Z])-\widehat{\mu}([Z]).
$$
The permutation we apply to the zigzags in homology class $[Z]$ is a so-called shuffle.
This means that $[Z]$ is split into disjoint intervals which are permuted,
while inside each interval the ordering is unchanged. The shuffle we apply to
the zigzags in $[Z]$ is depicted in Figure \ref{fig:shuffle};
the numbers $\widehat{\gl}([Z])$ etc. indicate the length of the interval.

Such shuffles must be applied to each homology class in the pattern
$\cZ=(Z_1,\ldots,Z_p)$. The composite result is the shuffle permutation
\begin{equation}\label{eq:shuffle}
\gs:\{1,\ldots,p\}\longrightarrow \{1,\ldots,p\}\,.
\end{equation}

\

\ssnnl{Definition.}\label{merging}
We define the \textit{merging matrix} $M_{\cA\cZ}$ for the set $\cA$ and the pattern of zigzags $\cZ$ to be the $p\times (p-q)$-matrix with $(i,j)$-entry
\begin{equation}\label{eq:merging matrix}
(M_{\cA\cZ})_{ij}=1\;\textrm{if}\;j=\gf(\gs(i))\,,\quad
(M_{\cA\cZ})_{ij}=0\;\textrm{if}\;j\neq\gf(\gs(i))\,,
\end{equation} 
with $\gf$ and $\gs$ as in (\ref{eq:pre-merge}) and (\ref{eq:shuffle}).

The \textit{merging step} in the algorithm multiplies the matrices
$B_\cZ$, $I_\cZ$, $P_\cZ$ and $Q_\cZ$ from the right with the matrix
$M_{\cA\cZ}$ and subsequently deletes 
the rows which correspond to intersection points in the pattern $\cZ$ which disappear in the merging process. These are recognized as the rows of $I_\cZ\,M_{\cA\cZ}$ with only one non-zero entry (namely $2$).

\textit{Thus the merging $\cZ\;\rightsquigarrow\;\cZ'$ is realized by}
\begin{equation}\label{eq:iteration 1}
\left\{
\begin{array}{rcl}
B_{\cZ'}\!\!\!\!&=&\!\!\!\!B_\cZ\,M_{\cA\cZ}
\\
I_{\cZ'}\!\!\!\!&=&\!\!\!\!\textrm{delete rows from }I_\cZ\,M_{\cA\cZ}
\\
P_{\cZ'}\!\!\!\!&=&\!\!\!\!\textrm{delete rows from }P_\cZ\,M_{\cA\cZ}
\\
Q_{\cZ'}\!\!\!\!&=&\!\!\!\!\textrm{delete rows from }Q_\cZ\,M_{\cA\cZ}
\end{array}\right.
\end{equation}

\

\ssnnl{}The transformation of patterns of zigzags $\cZ\;\rightsquigarrow\;\cZ'$
in (\ref{eq:iteration 1}) has been organized so that if $\cZ$ is a very good pattern, then
$\cZ'$ satisfies the Conditions \ref{def:zigzag pattern}.1-6,
\ref{def:good}.7 and the property formulated on the first line of Equation (\ref{eq:bands}). However,
$\cZ'$ need not satisfy \ref{def:good}.8, i.e. the equality
$$
|Z'_i\edot Z'_j|\,=\,\sharp(Z'_i\cap Z'_j)\,.
$$
need not hold for all pairs of zigzags $(Z'_i,Z'_j)$ in $\cZ'$.
This equality can only be violated if $Z'_i$ results from merging
zigzags $Z_{i_1}$ and $Z_{i_2}$ from $\cZ$ such that 
$Z_{i_1}\edot Z_{i_2}=1$ and
$(Z_{i_1}\edot Z'_j)(Z_{i_2}\edot Z'_j)<0$.
Since $[Z'_j]$ can only be the homology class of a zigzag in $\cZ$ or
the sum of two such, this can only happen if 
$[Z'_j]=\pm([Z_{i_1}]+[Z_{i_2}])=\pm[Z'_i]$.

The merging process in (\ref{eq:iteration 1}) is such that
if $Z'_i$ results from merging
zigzags $Z_{i_1}$ and $Z_{i_2}$ from $\cZ$, then every zigzag $Z'_j$ in the homology class $[Z'_i]$ (resp. in $-[Z'_i]$) is the result of merging 
a zigzag $Z_{j_1}$ from $[Z_{i_1}]$ (resp. $-[Z_{i_1}]$) with a zigzag $Z_{j_2}$ from $[Z_{i_2}]$ (resp. $-[Z_{i_2}]$).
Since $\cZ$ satisfies Condition \ref{def:good}.8. and 
$Z_{i_1}\edot Z_{i_2}=1$ we have in that case
$$
Z_{i_1}\cap Z_{j_1}=Z_{i_2}\cap Z_{j_2}=\emptyset\,,\qquad
\sharp(Z_{i_1}\cap Z_{j_2})=\sharp(Z_{i_2}\cap Z_{j_1})=1\,.
$$
Whence if $Z'_j\in\pm[Z'_i]$ and $Z'_i\neq Z'_j$, then
$\sharp(Z'_i\cap Z'_j)=2$.
\\
\textbf{Conclusion:}
\textit{For a pair of zigzags $(Z'_i,Z'_j)$ in $\cZ'$ we have:}
$$
|Z'_i\edot Z'_j|\,\neq\,\sharp(Z'_i\cap Z'_j)
\;\Longleftrightarrow\;
Z'_i\edot Z'_j=0\textrm { and } \sharp(Z'_i\cap Z'_j)=2\,.
$$

\

\ssnnl{}\label{clean}
We now transform the pattern of zigzags $\cZ'$ created in 
(\ref{eq:iteration 1}) into a very good one.
First we apply repairing moves 2 (see \ref{repairing 2}) to those pairs of 
zigzags $(Z'_i,Z'_j)$ which satisfy $[Z'_i]=[Z'_j]$ and
$|Z'_i\edot Z'_j|\,\neq\,\sharp(Z'_i\cap Z'_j)$ and which either both do or both do not
belong to a $+$-opposite pair (cf. Definition \ref{def:pair}).
Next we apply repairing moves 1 (see \ref{repairing 1}) wherever possible.
See also Example \ref{example:band} and Remark \ref{remark:reorder}
for the effect of repairing moves 2 and 1 on alternating sequences of opposite pairs. Thus with repairing moves 2 and 1 and possibly a permutation of columns in $I_{\cZ'}$, $P_{\cZ'}$, $Q_{\cZ'}$ we now have a pattern of zigzags which satisfies also the second half of Equation (\ref{eq:bands}) and in which two zigzags
$Z_i,\,Z_j$ with $[Z_i]=\pm[Z_j]$ do not intersect unless one is member of an 
opposite pair and the other is not.

Finally we apply repairing moves 3 (see \ref{repairing 3}) as follows for all homology classes $[Z]$ and $-[Z]$ which have resulted from merging and which satisfy $\sharp([Z])>\sharp(-[Z])>0$. In these circumstances the
zigzags in $-[Z]$ and the first $\sharp(-[Z])$ zigzags in $[Z]$
form a sequence $(Z_1,\ldots,Z_{2s})$ as in the beginning of \ref{repairing 3}.
For the zigzag $Z_0$ in \ref{repairing 3} we take the zigzag in $[Z]$ with the highest index. Note that while performing repairing move 3 we have first deleted from the matrices $I_\cZ$, $P_\cZ$ and $Q_\cZ$ all rows which corresponded with intersection points on one of the zigzags in the sequence $(Z_1,\ldots,Z_{2s})$.
Subsequently we have inserted rows for intersection points of a zigzag in the new sequence $(Z'_1,\ldots,Z'_{2s})$ with a zigzag which also intersects $Z'_0=Z_0$.
The previously applied repairing moves 2 had already removed all intersection points of $Z_0$ with the zigzags in $[Z]$ which were not in
$(Z_1,\ldots,Z_{2s})$.
So after applying repairing moves 3 two zigzags in $[Z]\cup(-[Z])$ do not intersect.

\

\ssnnl{} After the above merging and repairing moves we have produced a very good pattern of zigzags and now \textbf{return
to \ref{actual iteration} for the next iteration.}

\section{From $B_\cZ$, $I_\cZ$, $P_\cZ$, $Q_\cZ$ to $\KK_\cZ$ and back}\label{sec:KZ}
\ssnnl{}\label{QZAZ}
The conversion works for every good pattern $\cZ$ of, say $p$, zigzags.
Condition \ref{def:very good}.9 is not needed here.
The rows of $P_\cZ$ and $Q_\cZ$ must be taken modulo the row space of $B_\cZ$. This is achieved by multiplying $P_\cZ$ and $Q_\cZ$
from the right by a $p\times(p-2)$-matrix $A_\cZ$ with entries in $\ZZ$, such that $\rank(A_\cZ)=p-2$ and $B_\cZ\,A_\cZ=\nv$.

The rows of $P_\cZ\,A_\cZ$ represent points in $\modquot{\ZZ^p}{\ZZ^2B_\cZ}$.
We denote the set of these points by $\gB$, because these are in fact the black nodes in the dimer model. In the zigzag pattern these are the $+$-cells.
Similarly, we denote the set of rows of $Q_\cZ\,A_\cZ$ by $\gW$.
These are the white nodes in the dimer model and the $-$-cells in the zigzag pattern.

\

\ssnnl{Definition.}\label{def:kasteleyn}(cf. \cite{S1} Definition 8.2, \cite{S2} Definition 1)
The \textit{generalized Kasteleyn matrix} $\KK_\cZ(\vz,\vu)$ of a good pattern of zigzags $\cZ=(Z_1,\ldots,Z_p)$ is defined
as follows.
The rows of $\KK_\cZ(\vz,\vu)$ correspond $1:1$ with the elements of $\gB$
and the columns correspond $1:1$ with the elements of $\gW$.
The entries of $\KK_\cZ(\vz,\vu)$ are polynomials in two sets of variables
$\vz$ and $\vu$. The variables in $\vu=\{u_1,\ldots,u_p\}$ correspond $1:1$ 
with the zigzags in $\cZ$, and, hence, with the columns of $B_\cZ$, $I_\cZ$, 
$P_\cZ$ and $Q_\cZ$. The variables in $\vz=\{z_e\}$ correspond $1:1$ with the intersection points $e$ of $\cZ$ and, hence, with the rows of
$I_\cZ$, $P_\cZ$ and $Q_\cZ$.

For an intersection point $e$ we denote by $\vb(e)$ the element of $\gB$ which ``is'' the $e$-th row of $P_\cZ\,A_\cZ$. Similarly, $\vw(e)\in\gW$ ``is'' the $e$-th row of $Q_\cZ\,A_\cZ$. Finally, $i(e)$ and $j(e)$ are such that
$I_\cZ(e,i(e))=I_\cZ(e,j(e))=1$.

Finally, for $\vb\in\gB$ and $\vw\in\gW$ we define:
\begin{equation}\label{eq:kast}
\textit{the $(\vb,\vw)$-entry of
$\:\KK_\cZ(\vz,\vu)\:$ is}\;\sum_{e:\,\vb(e)=\vb,\,\vw(e)=\vw}z_e u_{i(e)}u_{j(e)}\,.
\end{equation}

\

\ssnnl{Example.} For the pattern of zigzags in Figure \ref{fig:zigzag matrices}
the generalized Kasteleyn matrix is:
$$
\KK_\cZ(\vz,\vu)\,=\,\left[
\begin{array}{ccc}
z_1u_1u_5&z_2u_1u_2+z_3u_4u_5&z_4u_2u_4\\
z_5u_1u_3&z_6u_1u_6+z_7u_3u_4&z_8u_4u_6\\
z_9u_3u_5&z_{10}u_5u_6+z_{11}u_2u_3&z_{12}u_2u_6
\end{array}
\right]\,.
$$

\

\ssnnl{Definition.}\label{def:complementary}
(cf. \cite{S2} Definition 3)
The \textit{complementary generalized Kasteleyn matrix} $\KK^c_\cZ(\vz,\vu)$ of a good pattern $\cZ$ of $p$ zigzags  is 
$$
\KK^c_\cZ(\vz,\vu)\,=\,
u_1\cdot\ldots\cdot u_p\KK^c_\cZ(\vz,u_1^{-1},\ldots,u_p^{-1})\,.
$$

\

\ssnnl{Remark.} The information in $\KK_\cZ(\vz,\vu)$ is in fact equivalent
with that in $B_\cZ$, $I_\cZ$, $P_\cZ$.
By Theorem 9.3 and \S 3.7 in \cite{S1} the columns of $B_\cZ$ are the primitive vectors along the sides of the Newton polygon of $\det\KK_\cZ(\vz,\vu)$
w.r.t. $u_1,\ldots,u_p$.
One recovers $P_\cZ$ and $I_\cZ$ as follows. Let $\cP$ and $\cQ$, respectively, be the sets of exponent vectors of the monomials in $u_1,\ldots,u_p$ which appear in the first columns of the matrices
$$
\KK_\cZ(\vz,\vu)\,\left(\KK^t_\cZ(\vz,\vu)\,\KK_\cZ(\vz,\vu)\right)^n
\quad\textrm{resp.}\quad 
\left(\KK^t_\cZ(\vz,\vu)\,\KK_\cZ(\vz,\vu)\right)^{n+1}
$$
for $n\in\ZZ_{\geq 0}$. Translations on $\ZZ^p$ by vectors in $\ZZ^2B_\cZ$
preserve $\cP$ and $\cQ$. Now take a `fundamental domain' 
$\cP^*\subset\cP$ for the $\ZZ^2B_\cZ$-action on $\cP$.
For each vector $\ga$ in $\cP^*$ let $\cQ_\ga$ be the set of vectors
$\gb$ in $\cQ$ such that $\ga-\gb$ has precisely two non-zero entries and these are both $1$. The rows of $P_\cZ$ resp. $I_\cZ$
are $\ga$ resp. $\ga-\gb$ with $\ga\in\cP^*$ and $\gb\in\cQ_\ga$.

\section{The principal $\cA$-determinant}\label{sec:Adet}

\ssnnl{The principal $\cA_\cZ$-determinant for a good pattern of zigzags.}\\
Let $\cZ=(Z_1,\ldots,Z_p)$ be a good pattern of zigzags and
let $\cA_\cZ$ denote the set of rows of the matrix $A_\cZ$ (see \ref{QZAZ}).
Let $\iota$ be the homomorphism
$$
\iota:\ZZ[z_e\,|\, e\textrm{ intersection point in } \cZ]\longrightarrow\ZZ\,,
\qquad \iota(z_e)=|Z_{i(e)}\edot Z_{j(e)}|\,.
$$
The \textit{principal $\cA_\cZ$-determinant} is defined in \cite{GKZ} and
Theorem 3 in \cite{S2} states that it is equal to
\begin{equation}\label{eq:AZdet}
\iota\left(\det \KK^c_\cZ(\vz,\vu)\right)\,.
\end{equation}
The relation
$$
I_\cZ\,=\,P_\cZ\,-\,Q_\cZ
$$
is precisely the one required in Condition 2 of \cite{S2}.

\

\ssnnl{Concluding remark about the principal $\cA$-determinant.}\\
For the principal $\cA$-determinant of
the set $\cA=\{\va_1,\ldots,\va_N\}$ in the Introduction we may have to make some slight adaptations to formula (\ref{eq:AZdet}), which reverse in a sense the transformation from $B_\cA$ to $B_\cZ$ in \ref{BA}. 

In (\ref{eq:AZdet}) the variables in $\vu=(u_1,\ldots,u_p)$ correspond $1:1$ with the columns of $B_\cZ$. Take a new set of variables 
$\vv=(v_1,\ldots,v_N)$ which correspond $1:1$ with the columns of $B_\cA$.
Recall that $B_\cZ$
is obtained from $B_\cA$ by permuting and splitting up columns as in (\ref{eq:primitive}). 
Then, to reverse the transformation from $B_\cA$ to $B_\cZ$
one must set $u_i=d_k\,v_k$ if the $i$-th column of $B_\cZ$ comes from the $k$-th column of $B_\cA$; here $d_k$ is the g.c.d. of the
two entries in the $k$-th column of $B_\cA$.

\end{document}